\documentclass[11pt]{article}   \usepackage{amsfonts}

\multiply\oddsidemargin by 2 \advance\oddsidemargin by -450bp
\advance\oddsidemargin by \textwidth \divide\oddsidemargin by 2

= 45bp

\begin{document}

\parindent = 0pt \baselineskip = 22pt   \parskip = \the\baselineskip

\newcommand{\Aa}{\mathbb{A}} \newcommand{\RR}{\mathbb{R}}
\newcommand{\FF}{\mathbb{F}} \newcommand{\CC}{\mathbb{C}}
\newcommand{\QQ}{\mathbb{Q}} \newcommand{\HH}{\mathbb{H}}
\newcommand{\ZZ}{\mathbb{Z}} \newcommand{\NN}{\mathbb{N}}

\newcommand{\Cc}{{\cal C}} \newcommand{\Ss}{{\cal S}}
\newcommand{\Real}{\mathop{\rm Re}} \newcommand{\Tr}{\mathop{\bf Tr}}
\newcommand{\llim}{\mathop{\rm l{.}i{.}m{.}}}

\newtheorem{thm}{Theorem}[section]

\newtheorem{lem}[thm]{Lemma}

\newtheorem{corollary}[thm]{Corollary}

\newtheorem{definition}[thm]{Definition}

\newtheorem{pre-note}[thm]{Note}
\newenvironment{note}{\begin{pre-note}\rm}{\end{pre-note}}

\newenvironment{proof}{\bf Proof\ \rm}{$\quad\clubsuit$}



\rightline{math/9904044}
\rightline{Includes math/9907103}
\rightline{11S80, 42B99, 47G10}
\vskip 2cm

\begin{center}
{\bf QUATERNIONIC GAMMA FUNCTIONS AND THEIR LOGARITHMIC DERIVATIVES  AS SPECTRAL
FUNCTIONS}\\

\vskip 0.5 true in

Jean-Fran\c{c}ois Burnol\\

November 2000\\

\end{center}

We establish Connes's local trace formula (related to the explicit formulae of
number theory) for the quaternions. This is done as an application of a study of
the central operator $H = \log(|x|) + \log(|y|)$ in the context of
invariant harmonic analysis. The multiplicative analysis of the additive Fourier 
transform gives a spectral interpretation to
generalized ``Tate Gamma functions'' (closely akin to the Godement-Jacquet
``$\gamma(s,\pi,\psi)$'' functions.) The analysis of $H$ leads furthermore to a spectral
interpretation for the logarithmic derivatives of these Gamma functions (which
are involved in ``explicit formulae''.)

\par

\vfill

Universit\'e de Nice\--\ Sophia Antipolis\\ Laboratoire J.-A. Dieudonn\'e\\ Parc
Valrose\\ F-06108 Nice C\'edex 02\\ France\\ burnol@math.unice.fr\\

\clearpage

\tableofcontents

\section{Introduction}

We establish for quaternions an analog of the trace formula obtained by Connes
in \cite{Co}  for a commutative local field $K$. This formula has the form
$\Tr(\widetilde{P_\Lambda}P_\Lambda\,U_f) =  2\log(\Lambda)f(1) + W(f) + o(1)$
(for $\Lambda\to\infty$), where $f$ is a test-function on $K^\times$, $U_f$ is
the operator of multiplicative convolution with $f$, $P_\Lambda$ and
$\widetilde{P_\Lambda}$ are cut-off projections (precise definitions will be
given later), all acting on the Hilbert space of square-integrable functions on
$K$. The contant term $W(f)$ was shown by Connes to be exactly the term arising
in the ``Weil's explicit formulae'' \cite{We1} of number theory.

We have shown in this abelian local case (see \cite{Bu1}, \cite{Bu2}, and the
related papers \cite{Bu3} and \cite{Bu4}) that the Weil term $W(f)$ can be
written as $-H(f)(1)$ for a certain dilaton-invariant operator $H$. We study
this operator in the non-commutative context of quaternions and then derive the
(analog of) Connes's asymptotic formula.   The proof would go through (with some
simplifications of course) equally well in the abelian case.

We first give some elementary lemmas of independent interest about self-adjoint
operators. We then study in multiplicative terms the additive Fourier Transform,
and this immediately leads to the definition of certain ``Quaternionic Tate
Gamma'' functions and to the analog of Tate's local functional equations
(\cite{Ta}, \cite{We2}). This is of course very much related to the
generalization to $GL(N)$ of Tate's Thesis in the work \cite{GJ} of
Godement-Jacquet (see \cite{Ja1}, \cite{Ja2} for reviews and references to
further works by other authors), where certain ``$\gamma(s,\pi,\psi)$''
functions, local $L$- and $\epsilon$-factors and associated functional equations
are studied. Also relevant is the classic monograph by Stein and Weiss \cite{St}
on harmonic analysis in euclidean spaces. In this paper we will follow a
completely explicit and accordingly elementary approach.

We introduce the ``conductor operator'' $H= \log(|x|) + \log(|y|)$ and show how
it gives an operator theoretic interpretation to the logarithmic derivatives of
the Gamma functions (which are involved in explicit formulae.) It is then a
simple matter to compute Connes's trace, and to obtain the asymptotic formula
$$\Tr(\widetilde{P_\Lambda}P_\Lambda\,U_f) = 2\log(\Lambda)f(1) - H(f)(1) +
o(1)$$ in a form directly involving our operator $H$. Further work leads to a
``Weil-like'' formulation for the constant term $H(f)(1)$, if so desired.

\section{Closed invariant operators}

It is well known that any bounded operator on $L^2(\RR,dx)$ which commutes with
translations is diagonalized by the additive Fourier transform (see for example
the Stein-Weiss monograph \cite{St}.) We need a generalization which applies to
(possibly) unbounded operators on Hilbert spaces of the form $L^2(G,dx)$ where
$G$ is a topological group. Various powerful statements are easily found in the
standard references on Hilbert spaces, usually in the language of spectral
representations of abelian von Neumann algebras. For lack of a reference
precisely suited to the exact formulation we will need, we provide here some
simple lemmas with their proofs.

\begin{definition} Let $L$ be a Hilbert space. A (possibly unbounded) operator 
$M$ on $L$ with domain $D$ is said to commute with the bounded operator $A$ if
$$\forall v\in L: v\in D\Rightarrow \Big(\ A(v)\in D\ \hbox{and}\ M(A(v)) =
A(M(v))\ \Big)$$
\end{definition}

\begin{thm}\label{L1}
Let $L$ be a Hilbert space and $G$ a (not necessarily abelian) group of unitary
operators on $L$. Let $\cal A$ be the von Neumann algebra of bounded operators
commuting with $G$. Let $M$ be a (possibly unbounded) operator on $L$, with
dense domain $D$. If the three following conditions are satisfied\\ (1)  $\cal
A$ is abelian\\ (2)  $(M,D)$ is symmetric\\ (3)  $(M,D)$ commutes with the
elements of $G$\\ then $(M,D)$ has a unique self-adjoint extension. This
extension commutes with the operators in $\cal A$.
\end{thm}


\begin{proof} We first replace $(M,D)$ by its double-adjoint so that we can
assume that $(M,D)$ is closed (it is easy to check that conditions (2) and (3)
remain valid). The problem is to show that it is self-adjoint. Let $K$ be the
range of the operator $M + i$. It is a closed subspace of $L$ (as $\|(M +
i)(\varphi)\|^2 = \|M(\varphi)\|^2 + \|\varphi\|^2$, and $M$ is closed). Let $R$
be the bounded operator onto $D$ which is orthogonal projection onto $K$
followed with the inverse of $M + i$. One checks easily that $R$ belongs to
$\cal A$, hence commutes with its adjoint $R^*$ which will also belong to $\cal
A$. Any vector $\psi$ in the kernel of $R$ is then in the kernel of $R^*$ (as
$<R^*\psi|R^*\psi>\ =\ <\psi|R\,R^*\psi>\ =\ 0$). So $\psi$ belongs to the
orthogonal complement to the range of $R$, that is $\psi = 0$ as the range of
$R$ is $D$. So $K = L$ and in the same manner $(M - i)(D) = L$. By the basic
criterion for self-adjointness (see \cite{Re}), $M$ is self-adjoint. Let
$A\in\cal A$. It commutes with the resolvant $R$ hence leaves stable its range
$D$. On $D$  one has $RA(M+i) = AR(M+i) = A = R(M+i)A$ hence $A(M+i) = (M+i)A$
so $A$ commutes with $M$.
\end{proof}

For the remainder of this section we let $G$ be a locally compact, Hausdorff,
topological \emph{abelian} group and $\widehat{G}$ its dual group. We refer to
\cite{Ru} for the basics of harmonic analysis on $G$. In particular we have a
Haar measure $dx$ (unique up to a multiplicative constant) and a Hilbert space
$L = L^2(G, dx)$. We also have a dual Haar measure $dy$ on $\widehat{G}$ such
that the Fourier transform $F(\varphi)(y) = \int \varphi(x) \overline{y(x)} dx$
is an isometry of $L$ onto $L^2(\widehat{G}, dy)$.  We sometimes identify the
two Hilbert spaces without making explicit the reference to $F$: so when we
write $f(y)\in L$ we really refer to $F^{-1}(f)\in L$ with $f\in
L^2(\widehat{G}, dy)$. No confusion should arise. We will assume that $dy$ is a
$\sigma$-finite measure so that there exists $\psi \in L$ with the property
$\psi(y)\neq 0\ a.e.$\ .

Let $a(y)$ be a measurable function on $\widehat{G}$, not necessarily
bounded. Let $D_a\subset L$ be the domain of square-integrable (equivalence
classes of measurable) functions $\varphi(x)$ on $G$ such that
$a(y)F(\varphi)(y)$ belongs to $L^2(\widehat{G}, dy)$. And let $M_a$ be the
operator with domain $D_a$ acting according to $\varphi\mapsto M_a(\varphi) =
F^{-1}(a\cdot F(\varphi))$. We write $a=b$ if the two functions $a(y)$ and
$b(y)$ are equal almost everywhere on $\widehat{G}$.

\begin{lem} The operator $(M_a, D_a)$ on $L^2(G,dx)$ commutes
with $G$. Furthermore $D_a$ is dense and $(M_a, D_a)$ is a closed operator. If
$(M_b, D_b)$ extends $(M_a,D_a)$, then in fact $a = b$ and $(M_b, D_b) =
(M_a,D_a)$. The adjoint of $(M_a,D_a)$ is $(M_{\overline{a}},D_{\overline{a}})$
(of course $D_{\overline{a}} = D_a$.)
\end{lem}

\begin{proof} We give the proof for completeness. The commutation with
  $G$-translations is clear. Then $D_a$  contains (the inverse
  Fourier transform of) $\psi(y)\over\sqrt{1 + |a(y)|^2}$ and all its
  translates. Hence if $f$ is orthogonal to $D_a$ then the function
  $\overline{f(y)}\psi(y)\over\sqrt{1 + |a(y)|^2}$ on $\widehat{G}$ belongs to
  $L^1(\widehat{G},dy)$ and has a vanishing ``inverse Fourier transform'', hence
  $f = 0$ (almost everywhere). It is also clear using $\psi(y)\over\sqrt{1 +
  |a(y)|^2}$ that if $(M_b, D_b)$ extends $(M_a,D_a)$, then $a = b$. Let us
  assume that the sequence $\varphi_j$ is such that $\varphi=\llim \varphi_j$
  and $\theta=\llim M_a(\varphi_j)$ both exist. Let us pick a pointwise on
  $\widehat{G}$ almost everywhere convergent subsequence
  $\varphi_{j_k}(y)$. Using Fatou's lemma we deduce that $\varphi$ belongs to
  $D_a$. Using Fatou's lemma again we get the vanishing of $\int_{\widehat{G}}
  |\theta(y) - a(y) \varphi(y)|^2\,dy$, and this shows that  $(M_a, D_a)$ is a
  closed operator. Finally let $f(y)$ be in the domain of the adjoint of
  $(M_a,D_a)$. There exists then an element $\theta$ of $L$ such that for any
  $\varphi \in D_a$ the equality $$\int f(y)\overline{a(y)\varphi(y)}\,dy = \int
  \theta(y)\overline{\varphi(y)}\,dy$$ holds. This implies that the two
  following functions of $L^1(\widehat{G}, dy)$:
$${f(y)\overline{a(y)\psi(y)}\over\sqrt{1 + |a(y)|^2}}\hbox{\quad
and\quad}{\theta(y)\overline{\psi(y)}\over\sqrt{1 + |a(y)|^2}}$$  have the same
Fourier transform on $G$, hence are equal almost everywhere. So $f \in
D_{\overline{a}}$ and $(M_a)^*(f) = (M_{\overline{a}})(f)$.
\end{proof}

\begin{thm}\label{T1} Let $(M, D)$ be a \emph{closed operator} on $L^2(G,dx)$
  commuting with $G$-translations. Then $(M, D) = (M_a, D_a)$ for a (unique)
  multiplicator $a$.
\end{thm}

\begin{note} For a bounded $M$ and $G = \RR$, this is proven in the classical
  monograph by Stein and Weiss \cite{St}, as a special case of a more general
  statement applying in $L^p$ spaces.
\end{note}

\begin{proof} 
Let us first assume that $M$ is bounded. We use the (inverse Fourier transform
of the) function $\psi(y)$ and define $a(y)$ to be
${M(\psi)(y)\over\psi(y)}$. Let us consider the domain $D$ consisting of all
finite linear combinations of translates of $\psi$. It is dense by the argument
using unicity of Fourier transform in $L^1$ we have used previously. Then $(M,
D) \subset (M_a, D_a)$, hence $(M_a, D_a)$ is also an extension of the closure
of $(M, D)$. As $M$ is assumed to be bounded this is $(M, L)$. But this means
that $D_a = L$ and that $M = M_a$ (we then note that necessarily $a$ is
essentially bounded).

The next case is when $M$ is assumed to be self-adjoint. Its resolvents $R_1 =
(M - i)^{-1}$ and $R_2 = (M + i)^{-1}$ are bounded and commute with $G$. Hence
they correspond to multiplicators $r_1(y)$ and $r_2(y)$. The kernel of $R_1$ is
orthogonal to the range of $R_2 = R_1^*$ which is all of $D$, so in fact it is
reduced to $\{0\}$. Hence $r_1(y)$ is almost everywhere non-vanishing. Let $f
\in D$ and $g = M(f)$. As $R_1(M(f) - i\,f) = f$ we get $g(y) = {1 +
i\,r_1(y)\over r_1(y)}\cdot f(y)$ and defining $a(y)$ to be ${1 + i\,r_1(y)\over
r_1(y)}$ we see that $(M_a, D_a)$ is an extension of $(M, D)$. Taking the
adjoints we deduce that $(M, D)$ is an extension of
$(M_{\overline{a}},D_{\overline{a}})$. So all three are equal (and $a$ is
real-valued).

For the general case we use the theorem of polar decomposition (see for example
\cite{Re}). There exists a non-negative self-adjoint operator $|M|$ with the
same domain as $M$ and a partial isometry $U$ such that $M = U|M|$. Further
conditions are satisfied which make $|M|$ and $U$ unique: so they also commute
with $G$. It follows from what was proven previously that $(M, D) \subset (M_a,
D_a)$ for an appropriate $a$ (the product of the multiplicators associated to
the self-adjoint $|M|$ and the bounded $U$). The adjoint $(M^*, D^*)$ also has a
dense domain and commutes with $G$, so in the same manner $(M^*, D^*) \subset
(M_b, D_b)$ for an appropriate $b$. The inclusion
$(M_{\overline{a}},D_{\overline{a}}) \subset (M^*, D^*) \subset (M_b, D_b)$
implies $b = \overline{a}$ and $(M_a, D_a) = (M, D)^{**}$. But the
double-adjoint coincides with the closed operator $(M,D)$.
\end{proof}

Let us mention an immediate corollary:

\begin{corollary} A closed symmetric operator on $L^2(G,dx)$ commuting with
$G$ is self-adjoint, and a symmetric operator which has a dense domain and
commutes with $G$ is essentially self-adjoint.
\end{corollary}

\section{Tate's functional equations}


Our first concern will be to introduce numerous notations. Let $\HH$ be the
space of quaternions with $\RR$-basis $\{1, i, j, k\}$ and table of
multiplication $i^2 = j^2 = k^2 = -1,\ ij = k = -ji,\ jk = i = -kj,\ ki = j =
-ik$. A typical quaternion will be denoted $x = x_0 + x_1 i + x_2 j + x_3 k$,
its conjugate $\overline{x} = x_0 - x_1 i - x_2 j - x_3 k$, its real part
$\Real(x) = x_0$, its (reduced) norm $n(x) = x\overline{x} = \overline{x}x =
x_0^2 + x_1^2 + x_2^2 + x_3^2$.

$\HH$ can also be considered as a left $\CC$-vector space with basis $\{1,
j\}$. We then write $a = x_0 + x_1 i$ and $b = x_2 + x_3 i$. Then $jaj^{-1} =
\overline{a}$, $x = a + bj$, and $n(x) = a\overline{a} + b\overline{b}$. The
action of $\HH$ on itself by right-multiplication sends $x = a + bj$ to the
$2\times2$ complex matrix
$$R_x = \pmatrix{a & -\overline{b} \cr b & \overline{a} \cr}$$  We write $V$ for
the complex vector space of complex-linear forms $\alpha: \HH\to\CC$. The forms
$A: x \mapsto a$ and  $B: x \mapsto b$ are a basis of $V$. We have a left action
of $\HH$ on $V$ with $x\in\HH$ acting as $\alpha(y) \mapsto \alpha(yx)$. This
left action represents the quaternion $x$ by the matrix
$$L_x = \pmatrix{a & b \cr -\overline{b} & \overline{a} \cr}$$ Also let $V_N =
\hbox{SYM}^N(V)$, for $N = 0, 1, \dots$ be the $N + 1$-dimensional complex
vector space with basis the monomials $A^j B^{N-j}, 0\leq j\leq N$.

Let $G = \HH^\times$ be the multiplicative group (with typical element $g$) and
$G_0 = \{ g \in G |\ n(g) = 1\}$ its maximal compact subgroup. Through the
assignment $g \mapsto L_g$ an isomorphism $G_0 \sim SU(2)$ is obtained, and the
$V_N$'s give the complete list of (isomorphism classes of) irreducible
representations of $G_0$.

The additive Fourier Transform $\cal F$ is taken with respect to the additive
character $x \mapsto \lambda(x) = e^{-2\pi\,i (x + \overline{x})}$. We note that
$\lambda(xy) = \lambda(yx)$. The choice we make for the normalization of $\cal
F$ is:
$${\cal F}(\varphi)(y) = \widetilde{\varphi}(y) = \int
\varphi(x)\lambda(-xy)\,dx$$ where $dx = 4dx_0dx_1dx_2dx_3$ is the unique
self-dual Haar measure for $\lambda$. With these choices the function $\omega(x)
= e^{-2\pi x\overline{x}}$ is its own Fourier transform.

\begin{definition}
The module $|g|$ of $g\in \HH^\times$ is defined by the equality of additive
Haar measures on $\HH$: $d(gx) = d(xg) = |g|dx$. It is expressed in terms of the
reduced norm by  $|g| = n(g)^2$.
\end{definition}

\begin{note}
The multiplicative (left- and right-) Haar measures on $G$ are the multiples of
$dg \over |g|$.
\end{note}

One has a direct product $G = (0, \infty) \times G_0$, $g = r g_0$, $r =
\sqrt{n(g)} = |g|^{1/4}$. We write $d\sigma$ for the Euclidean surface element
on $G_0$ (for the coordinates $x_i$), so that $dx = 4r^3\, drd\sigma$. The rule
for integrating functions of $r$ is $\int g(r) dx = 8\pi^2\,\int_0^\infty g(r)
r^3 dr$ as is checked with $\omega(x)$.   So $d\sigma = 2\pi^2\,d^*g_0$ where
$d^*g_0$ is the Haar measure on $G_0$ with total mass $1$.

\begin{definition}
The normalized Haar measure on $G$ is defined to be  $d^*g = {1 \over 2
\pi^2}{dg \over |g|} = 4 {dr\over r}\,d^*g_0$. It is chosen so that its
push-forward under the module map $g\mapsto u = |g|\in\RR^{\times+}$ is ${du
\over u} = 4 {dr\over r}$.
\end{definition}

The multiplicative group $G$ acts in various unitary ways on $L^2 := L^2(\HH,
dx)$:\\ \centerline{\hfill$L_1(g) : \varphi(x) \mapsto
|g|^{1/2}\varphi(xg)$\hfill $R_1(g) : \varphi(x) \mapsto
|g|^{1/2}\varphi(gx)$\hfill} and also  $L_2(g) = R_1(g^{-1})$ and $R_2(g) =
L_1(g^{-1})$.

\begin{definition}
The \emph{Inversion} $I$ is the unitary operator on $L^2(\HH,dx)$ acting as
$\varphi(x) \mapsto {1\over|x|}\varphi({1\over x})$. The \emph{Gamma operator}
is the composite $\Gamma={\cal F}I$.
\end{definition}

\begin{thm}
The Gamma operator commutes with both left actions $L_1$ and $L_2$ and with both
right actions $R_1$ and $R_2$ of $G$ on $L^2$.
\end{thm}

\begin{proof}
One just checks that ${\cal F}$ intertwines $L_1$ with $L_2$, and also $R_1$
with $R_2$ and that the inversion $I$ also intertwines $L_1$ with $L_2$, and
$R_1$ with $R_2$.
\end{proof}

\begin{definition}
The \emph{basic isometry} is the map $\phi(x) \mapsto f(g) = \sqrt{2\pi^2\,|g|}
\; \phi(g)$ between $L^2(\HH,dx)$ and $L^2(G, d^*g)$.
\end{definition}

\begin{note}
It is convenient to avoid using any notation at all for the basic isometry. So
we still denote by ${\cal F}$ the additive Fourier transform transported to the
multiplicative setting. The inversion $I$ becomes $f(g) \mapsto f(g^{-1})$. The
Gamma operator is still denoted $\Gamma$ when viewed as acting on $L^2(G, d^*g)$.
\end{note}

The spectral decomposition of $L^2((0, \infty),{du\over u})$ is standard Fourier
(or Mellin) theory (alternatively we can apply Theorem {\bf\ref{T1}} here): any bounded
operator $M$ commuting with multiplicative translations is given by a measurable
bounded multiplier $a(\tau)$ in dual space $L^2(\RR, {d\tau \over 2\pi})$:
$$G_1(u) = \llim_{\Lambda \rightarrow \infty} \int_{-\Lambda}^{\Lambda}
\psi(\tau) u^{-i\tau} {d\tau \over 2\pi} \Longrightarrow M(G_1)(u) =
\llim_{\Lambda \rightarrow \infty} \int_{-\Lambda}^{\Lambda} a(\tau)\psi(\tau)
u^{-i\tau} {d\tau \over 2\pi}$$

 On the other hand the spectral decomposition of $L^2(G_0,\,d^*g_0)$ is part of
 the Peter-Weyl theory: it tells us that $L^2(G_0,\,d^*g_0)$ decomposes under
 the $L_1 \times R_1$ action by $G_0 \times G_0$ into a countable direct sum
 $\oplus_{N\geq0} W_N$ of finite dimensional irreducible, non-isomorphic,
 modules. This is also the isotypical decomposition under either $L_1$ alone or
 $R_1$ alone (for which $W_N$ then contains $N+1$ copies of $V_N$.) Using the
 standard theory of tensor products of separable Hilbert spaces (see for example
 \cite{Re})  we have:

\begin{lem} The isotypical
decomposition of $L^2(G,\,d^*g)$ under the compact group $G_0 \times G_0$ acting
through $L_1 \times R_1$ is
$$L^2(G,\,d^*g) = L^2((0, \infty),{du\over u}) \otimes L^2(G_0,\,d^*g_0) =
\oplus_N L^2((0, \infty),{du\over u})\otimes W_N$$
\end{lem}

\begin{lem}\label{L2}
Let $M$ be a bounded operator on $L^2$ which commutes with both the $L_1$ and
$R_1$ actions of $G$. Then to each integer $N\geq 0$ is associated an
(essentially bounded) multiplicator $a_N(\tau)$ on $\RR$, unique up to equality
almost everywhere, such that
$$\psi\in L^2(\RR, {d\tau \over 2\pi}), \  G_1(u) = \llim_{\Lambda \rightarrow
\infty} \int_{-\Lambda}^{\Lambda} \psi(\tau) u^{-i\tau} {d\tau \over 2\pi}$$
$$\Rightarrow \forall F\in W_N\quad M(FG_1) = FG_2$$
$$\hbox{with }G_2(u) = \llim_{\Lambda \rightarrow \infty}
\int_{-\Lambda}^{\Lambda} a_N(\tau)\psi(\tau) u^{-i\tau} {d\tau \over 2\pi}$$
and where $FG_1$ is the function $g\mapsto F({g\over|g|^{1/4}})G_1(|g|)$ and
$FG_2$ the function $g\mapsto F({g\over|g|^{1/4}})G_2(|g|)$.
\end{lem}

\begin{proof}
Let us take $f \in L^2((0, \infty),{du\over u})$ and consider the linear
operator on $L^2(G_0,d^*g_0)$:
$$F(g_0) \mapsto \left(g_0 \mapsto \int_0^\infty \overline{f(u)}M(f \otimes
F)(g_0\,u^{1/4}){du\over u}\right)$$ It commutes with the action of $G_0 \times
G_0$ hence stabilizes each $W_N$ and is a multiple $a_N^f$ of the identity
there. On the other hand, if we choose $F_1$ and $F_2$ in $W_N$ and consider
$$f \mapsto \left(u \mapsto \int_{G_0} \overline{F_2(g_0)} M(f \otimes
F_1)(g_0\,u^{1/4})\,d^*g_0\right)$$ we obtain a bounded operator $M(F_1,F_2)$ on
$L^2((0, \infty),{du\over u})$ commuting with dilations and such that
$$<f |M(F_1,F_2)(f)> = <F_2 | M_N^f(F_1) > = a_N^f <F_2 | F_1>$$ where the
let-hand bracket is computed in $L^2((0, \infty),{du\over u})$ while the next
two are  in $L^2(G_0,\,d^*g_0)$. So $M(F_1,F_2)$ depends on $(F_1,F_2)$ only
through $<F_2 | F_1>$. We then let $a_N(\tau)$ be the spectral multiplier
associated to $M(F, F)$ for an arbitrary $F$ satisfying $<F|F> = 1$.
\end{proof}

\begin{corollary}
The von Neumann algebra $\cal A$ of bounded operators commuting simultaneously
with the left and right actions of the multiplicative quaternions on $L^2(\HH,
dx)$ is abelian.
\end{corollary}

\begin{lem}
A self-adjoint operator $M$ commuting with both left and right actions of $G$
commutes with any operator of the von Neumann algebra $\cal A$. In particular it
commutes with $\Gamma$.
\end{lem}

\begin{proof}
One applies Theorem {\bf\ref{L1}}.
\end{proof}

\begin{definition}\label{D1}
The \emph{quaternionic Tate Gamma functions} are the multiplicators
$\gamma_N(\tau)$ ($N\geq 0$) associated to the unitary operator $\Gamma$.
\end{definition}

\begin{note}
This generalizes the Gamma functions of Tate for $K=\RR$ anf $K=\CC$
(\cite{Ta}). In all cases they are indexed by the characters of the maximal
compact subgroup of the multiplicative group $K^\times$.
\end{note}

\begin{lem}
There is a smooth function in the equivalence class of
$\gamma_N(\tau)$.
\end{lem}

\begin{proof}
If the function $G_1(u)$ on $(0, \infty)$ is chosen smooth with compact support
(so that $\psi(\tau)$ is entire) then, for any $F\in W_N$ the function $FG_1$,
viewed in the additive picture, is smooth on $\HH$, has compact support, and
vanishes identically in a neighborhood of the origin. So its image under the
inversion also belongs to the Schwartz class in the additive picture on
$\HH$. Hence $\Gamma(FG_1)$ can be written as  $|g|^{1/2} \phi(g)$ for some
Schwartz function $\phi(x)$ of the additive variable $x$. One checks that this
then implies that $G_2(u)$ is a Schwartz function of the variable $\log(u)$ (we
assume that $F$ does not identically vanish of course), hence that
$\gamma_N(\tau) \psi(\tau)$ is a Schwartz function of $\tau$. The various
allowable $\psi$'s have no common zeros so the conclusion follows.
\end{proof}

\begin{note}
From now on $\gamma_N$ refers to this unique smooth representative. It is
everywhere of modulus $1$ as $\Gamma$ is a unitary operator.
\end{note}

\begin{note}
Any function  $F \in W_N$ will now be considered as a function on all of $G =
\HH^\times$ after extending it to be constant along each radial line. It is not
defined at $x=0$ of course.
\end{note}

Let $F\in W_N$. For $\Real(s) > 0$, $F(x) |x|^{s-1}$ is a tempered distribution
on $\HH$, hence has a distribution-theoretic Fourier Transform. At first we only
consider $s = {1\over 2} + i\tau$:

\begin{lem}\label{L4}
As distributions on $\HH$
$${\cal F}(F({1\over x})\,|x|^{-{1\over 2} + i\tau}) = \gamma_N(\tau)
F(x)\,|x|^{-{1\over 2} - i\tau}$$
\end{lem}

\begin{proof}
 We have to check the identity:
$$\int F({1\over y})\,|y|^{-{1\over 2} + i\tau}\widetilde{\varphi}(y)\,dy =
\gamma_N(\tau)\cdot\int F(x)\,|x|^{-{1\over 2} - i\tau}\varphi(x)\,dx$$ for all
Schwartz functions $\varphi(x)$ with Fourier Transform
$\widetilde{\varphi}(y)$. Both integrals are analytic in $\tau\in\RR$, hence
both sides are smooth (bounded) functions of $\tau$. It will be enough to prove
the identity after integrating against $\psi(\tau)\, {d\tau \over 2\pi}$ with an
arbitrary Schwartz function $\psi(\tau)$. With the notations of Lemma
{\bf\ref{L2}}, we have to check
$$\int F({1\over y})G_1({1\over y})|y|^{- {1\over 2}}\widetilde{\varphi}(y)\,dy
=  \int F(x)\,G_2(x)|x|^{-{1\over 2}}\varphi(x)\,dx$$ But, by Lemma
{\bf\ref{L2}}, and by Definition {\bf\ref{D1}}, $F(x)\,G_2(x)|x|^{-{1\over 2}}$
is just the Fourier Transform in $L^2(\HH, dx)$ of $F({1\over y})G_1({1\over
y})|y|^{- {1\over 2}}$, so this reduces to the $L^2$-identity
$$\int \psi(y)\widetilde{\varphi}(y)\,dy = \int
\widetilde{\psi}(x)\varphi(x)\,dx$$
\end{proof}

\begin{thm}\label{T2}
Let $F\in W_N$. There exists an analytic function $\Gamma_N(s)$ in $0 < \Real(s)
<1$ depending only on $N\in \NN$ and such that the following identity of
tempered distributions on $\HH$ holds for each $s$ in the critical strip
($0<\Real(s)<1$):
$${\cal F}(F({1\over x})\,|x|^{s -1}) = \Gamma_N(s) F(x)\,|x|^{-s}$$
\end{thm}

\begin{proof}
We have to check an identity:
$$\int F({1\over y})\,|y|^{s-1}\widetilde{\varphi}(y)\,dy = \Gamma_N(s)\cdot\int
F(x)\,|x|^{-s}\varphi(x)\,dx$$ for all Schwartz functions $\varphi(x)$ with
Fourier Transform $\widetilde{\varphi}(y)$. Both integrals are analytic in the
strip $0 < \Real(s) <1$, their ratio is thus a meromorphic function, which
depends neither on $F$ nor on $\varphi$ as it equals $\gamma_N(\tau)$ on the
critical line. Furthermore for any given $s$ we can choose $\varphi(x) =
\overline{F(x)} \alpha(|x|)$, with $\alpha$ having very small support around
$|x| = 1$ to see that this ratio is in fact analytic.
\end{proof}

\begin{note}
This is the analog for quaternions of Tate's ``local functional equation''
\cite{Ta}, in the distribution theoretic flavor advocated by Weil \cite{We2}. We
followed a different approach than Tate, as his proof does not go through that
easily in the non-commutative case.
\end{note}

Let $\Gamma(s)$ be Euler's Gamma function ($\int_0^\infty e^{-u}u^s\,{du\over
u}$).

\begin{thm}\label{T3} We have for each $N\in\NN$:
$$\Gamma_N(s) = i^N (2\pi)^{2-4s} {\Gamma(2s + {N\over 2})\over\Gamma(2(1-s) +
{N\over 2})}$$
\end{thm}

\begin{proof}
Let $0\leq j \leq N$ and $\omega_j(x) = \overline{A(x)}^{N-j}\overline{B(x)}^j
e^{-2\pi x\overline{x}} = \overline{a}^{N-j}\,\overline{b}^j\,\omega(x)$. One
checks that $\widetilde{\omega_j}(y) =
(-1)^j\,i^N\,{\alpha}^{N-j}\,\overline{\beta}^j \, \omega(y)$ ($y = \alpha +
\beta j$). We choose as homogeneous function $F_j(x) =
a^{N-j}\,b^j\,|x|^{-N/4}$. For these choices the identity of Theorem
{\bf\ref{T2}} becomes
$$i^N \int (\alpha\overline{\alpha})^{N-j}(\beta\overline{\beta})^{j} e^{-2\pi
y\overline{y}} |y|^{s - 1 - N/4} dy = \Gamma_N(s) \int
(a\overline{a})^{N-j}(b\overline{b})^{j} e^{-2\pi x\overline{x}} |x|^{-s - N/4}
dx$$  Adding a suitable linear combinations of these identities for $0\leq j
\leq N$ gives
$$i^N \int (y\overline{y})^{N}\,e^{-2\pi y\overline{y}} |y|^{s - 1 - N/4} dy =
\Gamma_N(s) \int (x\overline{x})^{N}\,e^{-2\pi x\overline{x}} |x|^{-s - N/4}
dx$$ hence the result after evaluating the integrals in terms of $\Gamma(s)$.
\end{proof}

\section{The central operator $H = \log(|x|) + \log(|y|)$}

\begin{definition}
We let $\Delta \subset {\cal C}^\infty(G)$ be the vector space of finite linear
combinations of functions $f(g) = F(g_0)K(\log(|g|))$ with $F$ in one of the
$W_N$'s (hence smooth) and $K$ a Schwartz function on $\RR$. It is a dense
sub-domain of $L^2$.
\end{definition}

\begin{thm} $\Delta$ is stable under ${\cal F}$.\end{thm}

\begin{proof}
We have to show that $\gamma_N(\tau)$ is a multiplier of the Schwartz class. Let
$h_N(\tau) = -i {\gamma_N^\prime(\tau) \over \gamma_N(\tau)}$. Using Theorem {\bf\ref{T3}}
and the partial fraction expansion of the logarithmic derivative of $\Gamma(s)$
(as in \cite{Bu1} for the real and complex Tate Gamma functions), or
Stirling's formula, or any other means, one finds $h_N(\tau) =
O(\log(1+|\tau|))$, $h_N^{(k)}(\tau) = O(1)$, so that $\gamma_N^{(k)}(\tau) =
O(\log(1+|\tau|)^k)$.
\end{proof}

Let $A$ be the operator on $L^2(\HH,dx)$ of multiplication with $\log(|x|)$. As
it is unbounded, we need a domain and we choose it to be $\Delta$. Of course
$(A, \Delta)$ is essentially self-adjoint. It is unitarily equivalent to the
operator $(B, \Delta)$, $B = {\cal F}A{\cal F}^{-1}$. Clearly:

\begin{lem} The domain $\Delta$ is stable under $A$ and $B$.\end{lem}

\begin{definition}
The \emph{conductor operator} is the operator $H = A + B$:
$$H = \log(|x|) + \log(|y|)$$ This is  an unbounded operator defined initially
on the domain $\Delta$.

\end{definition}

\begin{lem} The conductor operator $(H,\Delta)$ commutes with the left and with
  the right actions of $G$ and is symmetric. 
\end{lem}

This is clear. Applying now Theorem {\bf\ref{L1}} we deduce:

\begin{thm} The conductor operator $(H,\Delta)$ has a unique self-adjoint extension.
\end{thm}

We will simply denote by $H$ and call ``conductor operator'' this self-adjoint
extension.

\begin{thm} The conductor operator $H$ commutes with the inversion $I$.
\end{thm}

\begin{proof} Indeed, it commutes with $\Gamma$ by Theorem {\bf\ref{L1}} and it
  commutes with $\cal F$ by construction. 
\end{proof}

We now want to give a concrete description of its spectral functions.

\begin{definition} Let for each $N\in\NN$ and $\tau\in\RR$:
$$h_N(\tau) = -i {\gamma_N^\prime(\tau) \over \gamma_N(\tau)}$$
$$k_N(\tau) = - h_N^\prime(\tau)$$
\end{definition}

Explicit computations prove that the functions $h_N$ are left-bounded ($\exists
C\,\forall\tau\,\forall N\ h_N(\tau)\geq -C$) and that the functions $k_N$ are
bounded ($\exists C\,\forall\tau\,\forall N\ |k_N(\tau)|\leq C$.) We need not
reproduce these computations here, which use only the partial fraction expansion
of Euler's gamma function, as similar results are provided in \cite{Bu1} in the
real and complex cases.

Let $f(g) = F(g_0)\phi(|g|)$ be an element of $\Delta$, $F \in W_N \subset
L^2(G_0, dg_0)$, $\phi\in L^2((0,\infty), {du\over u})$, $\phi$ being a Schwartz
function of $\log(u)$ ($u = |g|$). We can also consider $f$ to be given as a
pair $\{F, \psi\}$ with $\psi(\tau) = \int_0^\infty \phi(u) u^{i\tau} {du\over
u}$ being a Schwartz function of $\tau$. Then $A(f)$ is given by the pair $\{F,
D(\psi)\}$ where $D$ is the differential operator ${1\over i}{d\over
d\tau}$. This implies that $\Gamma A\Gamma^{-1} (f)$ corresponds to the pair
$\{F, D(\psi) - h_N\cdot \psi\}$. On the other hand  $\Gamma A\Gamma^{-1} = - B$
so $H(f)$ corresponds to the pair $\{F, h_N\cdot \psi\}$. The commutation with
the inversion $I$ translates into $h_N(-\tau) = h_N(\tau)$. Also: $K = i[B, A] =
-i[A, H]$ sends the pair $\{F, \psi\}$ to $\{F, k_N\cdot \psi\}$, hence is
bounded and anti-commutes with the Inversion. We have proved:

\begin{thm}\label{T4}
The operator $\log(|x|) + \log(|y|)$ is  self-adjoint, left-bounded, commutes
with the left- and right- dilations, commutes with the Inversion, and its
spectral functions are the functions $h_N(\tau)$. The operator $i\,[\log(|y|),
\log(|x|)]$ is bounded, self-adjoint, commutes with the left- and right-
dilations, anti-commutes with the Inversion and its spectral functions are the
functions $k_N(\tau)$.
\end{thm}

We now conclude this chapter with a study of some elementary
distribution-theoretic properties of $H$. For this we need the analytic
functions of $s$ ($0 < \Real(s) < 1$) indexed by $N\in\NN$:
$$H_N(s) = {d\over ds} \log(\Gamma_N(s))$$ (so that $h_N(\tau) = H_N({1\over 2}
+ i\tau)$).

\begin{lem}
Let $\varphi(x)$ be a Schwartz function on $\HH$. Then $H(\varphi)$ is
continuous on $\HH \backslash\{0\}$, is $O(\log(|x|))$ for $x \rightarrow 0$,
and is $O(1/|x|)$ for $|x| \rightarrow\infty$. Furthermore, for any $F \in W_N$
(constant along radial lines), the following identity holds for $0 < \Real(s) <
1$:
$$\int H(\varphi)(x) F(x) |x|^{-s} dx = H_N(s) \int \varphi(x) F(x) |x|^{-s} dx$$
\end{lem}

\begin{proof}
Assuming the validity of the estimates we see that both sides of the identity
are analytic functions of $s$, so it is enough to prove the identity on the
critical line:
$$\int H(\varphi)(x) F(x) |x|^{-{1\over 2} - i\tau} dx = h_N(\tau) \int
\varphi(x) F(x) |x|^{-{1\over 2} - i\tau} dx$$ 
As in the proof of Lemma {\bf\ref{L4}}, it
is enough to prove it after integrating against an arbitrary Schwartz function
$\psi(\tau)$. With $G(u) = \int_{-\infty}^\infty \psi(\tau) u^{-i\tau} {d\tau
\over 2\pi}$, and using Theorem {\bf\ref{T4}} this becomes
$$\int H(\varphi)(x) F(x) G(|x|) |x|^{-{1\over 2}} dx = \int \varphi(x) H(FG)(x)
|x|^{-{1\over 2}} dx$$ (on the right-hand-side $H(FG)$ is computed in the
multiplicative picture, on the left-hand-side $H(\varphi)$ is evaluated in the
additive picture). The self-adjointness of $H$ reduces this to
$\overline{H(\overline{\varphi})} = H(\varphi)$, which is a valid identity.

For the proof of the estimates we observe that $B(\varphi)$ is the Fourier
transform of an $L^1$-function hence is continuous, so that we only need to show
that it is $O(1/|x|)$ for $|x| \rightarrow \infty$. For this we use that
$B(\varphi)$ is additive convolution of $-\varphi$ with the distribution $G =
{\cal F}(-\log(|y|))$. The estimate then follows from the formula for $G$ given
in the next lemma.
\end{proof}

\begin{lem}\label{L5} The distribution $G(x) = {\cal F}(-\log(|y|))$ is given as:
$$G(\varphi) = \int_{|x| \leq 1}(\varphi(x) - \varphi(0))\, {dx \over
2\pi^2\,|x|} + \int_{|x| > 1}\varphi(x)\, {dx \over 2\pi^2\,|x|}+\ (4\log(2\pi)
+ 4\gamma_e - 2)\varphi(0)$$
\end{lem}

\begin{proof}
We have used the notation $\gamma_e = - \Gamma^\prime(1)$ for the
Euler-Mascheroni's constant ($=0.577\dots$). Let $\Delta_s$ for $\Real(s)>0$ be
the homogeneous distribution $|x|^{s-1}$ on $\HH$. It is a tempered
distribution. The formula
$$\Delta_s(\varphi) = \int_{|x| \leq 1}(\varphi(x) - \varphi(0))|x|^{s-1}\, dx +
\int_{|x| > 1}\varphi(x)|x|^{s-1}\, dx + {2\pi^2\over s}\varphi(0)$$ defines its
analytic continuation to $\Real(s) > -{1 \over 4}$, with a simple pole at $s =
0$. Using
$${\cal F}(\Delta_s) = \Gamma_0(s) \Delta_{1-s}$$ for $s = 1 - \varepsilon$,
$\varepsilon \rightarrow 0$, and expanding in $\varepsilon$ gives
$$\varphi(0) + \varepsilon G(\varphi) + O(\varepsilon^2) = \Gamma_0(1 -
\varepsilon)\cdot\left\{{2\pi^2\over \varepsilon}\varphi(0) + \int_{|x| \leq
1}(\varphi(x) - \varphi(0))\, {dx\over|x|} + \int_{|x| > 1}\varphi(x)\,
{dx\over|x|} + O(\varepsilon) \right\}$$ As $\Gamma_0(1 - \varepsilon) = {1\over
2\pi^2}(\varepsilon + (4\log(2\pi) + 4\gamma_e - 2) \varepsilon^2 +
O(\varepsilon^3))$ the result follows.
\end{proof}

\section{The trace of Connes for quaternions}

Let $f(g)$ be a smooth function with compact support on $\HH^\times$. Let $U_f$
be the bounded operator $\int f(g) L_2(g)\,d^*g$ on $L^2(\HH, dx)$ of left
multiplicative convolution. So
$$U_f: \varphi(x) \mapsto \int_G f(g){1\over\sqrt{|g|}}\varphi(g^{-1}x)\,d^*g$$
The composition $U_f\,{\cal F}$ of $U_f$ with the Fourier Transform ${\cal F}$
acts as
\begin{eqnarray*}
\varphi(x) &\mapsto& \int_G\int_\HH
f(g){1\over\sqrt{|g|}}\lambda(-g^{-1}xy)\varphi(y)\,dy\,d^*g \\
&=& \int_\HH\int_G
f({1\over
g}){\sqrt{|g|}}\lambda(-gxy)\,d^*g\,\varphi(y)\,dy\\ 
&=&
{1\over\sqrt{2\pi^2}}\int_{y\in\HH} \left( \int_{Y\in\HH} f({1\over
Y}){1\over\sqrt{2\pi^2|Y|}}\lambda(-Yxy)\,dY\right) \varphi(y)\,dy \\ 
&=&
{1\over\sqrt{2\pi^2}}\int_\HH {\cal F}(I(f)_a)(xy)\varphi(y)\,dy\\
\end{eqnarray*}
In this last equation $I(f)_a$ is the additive representative
${1\over\sqrt{2\pi^2|Y|}}f({1\over Y})$ of $I(f)$. Finally denoting similarly
with $\Gamma(f)_a$ the additive representative of $\Gamma(f)$ we obtain
$$(U_f{\cal F})(\varphi)(x) = {1\over\sqrt{2\pi^2}}\int_\HH {\Gamma(f)_a}(xy)\varphi(y)\,dy$$

As $f$ has compact support on $\HH^\times$ we note that $I(f)_a$ is smooth with
compact support on $\HH$ and that $\Gamma(f)_a$ belongs to the Schwartz
class. Following Connes (\cite{Co}, for $\RR$ or $\CC$ instead of $\HH$), our
goal is to compute the trace $\Tr(\Lambda)$ of the operator
$\widetilde{P_\Lambda}P_\Lambda\,U_f$, where $\widetilde{P_\Lambda} = {\cal
F}P_\Lambda{\cal F}^{-1}$ and $P_\Lambda$ is the cut-off projection to functions
with support in $|x| \leq \Lambda$. Our reference for trace-class operators will
be \cite{Go}. We recall that if $A$ is trace-class then for any bounded $B$,
$AB$ and $BA$ are trace-class and have the same trace. Also if $K_1$ and $K_2$
are two Hilbert-Schmidt operators given for example as $L^2-$kernels $k_1(x,y)$
and $k_2(x,y)$ on a measure space $(X, dx)$ then $A = K_1^*\, K_2$ is
trace-class and its trace is the Hilbert-Schmidt scalar products of $K_1$ and
$K_2$:
$$\Tr(K_1^*\, K_2) = \int\int \overline{k_1(x,y)}\,k_2(x,y)\ dxdy$$

The operator $P_\Lambda {\cal F}^{-1} P_\Lambda$ is an operator with kernel a
smooth function restricted to a finite box (precisely it is $\lambda(xy),\
|x|,|y|\leq\Lambda$). Such an operator is trace class, as is well-known (one
classical line of reasoning is as follows: taking a smooth function $\rho(x)$
with compact support, identically $1$ on $|x|\leq\Lambda$, and $Q_\rho$ the
multiplication operator with $\rho$, one has $P_\Lambda {\cal F}^{-1} P_\Lambda
= P_\Lambda Q_\rho {\cal F}^{-1} Q_\rho P_\Lambda$, so that it is enough to
prove that $Q_\rho {\cal F}^{-1} Q_\rho$ is trace-class. This operator has a
smooth kernel with compact support, so we can put the system in a box, and
reduce to an operator $K$ with smooth kernel on a torus. Then $K = (1 +
\Delta)^{-n}(1 + \Delta)^{n}K$ with $\Delta$ the positive Laplacian. For $n$
large enough, $(1 + \Delta)^{-n}$ is trace-class, while $(1 + \Delta)^{n}K$ is
at any rate bounded.)\par

So Connes's operator $\widetilde{P_\Lambda}P_\Lambda\,U_f = {\cal F}\cdot
P_\Lambda {\cal F}^{-1} P_\Lambda\cdot U_f$ is indeed trace class and
$$\Tr(\widetilde{P_\Lambda}P_\Lambda\,U_f) = \Tr(P_\Lambda {\cal F}^{-1}
P_\Lambda\cdot U_f{\cal F}) = \Tr(P_\Lambda {\cal F}^{-1} P_\Lambda\cdot
P_\Lambda U_f{\cal F} P_\Lambda)$$ can be computed as a Hilbert-Schmidt scalar
product:
$$\Tr(\widetilde{P_\Lambda}P_\Lambda\,U_f) = {1\over\sqrt{2\pi^2}}\int\int_{|x|,
|y| \leq\Lambda}\lambda(xy)\Gamma(f)_a(xy)\,dxdy$$ using the change of
variable $(x,y) \mapsto (Y= xy, y)$
$$\Tr(\widetilde{P_\Lambda}P_\Lambda\,U_f) = \sqrt{2\pi^2}
\int_{|Y|\leq\Lambda^2} \lambda(Y) \Gamma(f)_a(Y)\left(\int_{{|Y|\over
\Lambda}\leq|y|\leq\Lambda} {dy\over 2\pi^2 |y|}\right)\,dY$$
$$\Tr(\widetilde{P_\Lambda}P_\Lambda\,U_f)  = \sqrt{2\pi^2}
\int_{|Y|\leq\Lambda^2} \Big(2\log(\Lambda) - \log(|Y|)\Big)\lambda(Y)
\Gamma(f)_a(Y)\,dY\leqno\bf(C)$$  
This integral is an inverse (additive) Fourier transform evaluated at $1$. As
$\Gamma = {\cal F}I$ itself involves a Fourier transform the final result is
just  $\sqrt{2\pi^2}M_\Lambda(I(f)_a)(1)$ where $M_\Lambda$ is the self-adjoint
operator $(2\log(\Lambda) - B)_+ = \max(2\log(\Lambda) - B, 0)$. If we recall
that $\sqrt{2\pi^2}$ is involved in the basic isometry from the additive to the
multiplicative picture, we can finally express everything back in the
multiplicative picture:

\begin{thm} The Connes operator $\widetilde{P_\Lambda}P_\Lambda\,U_f$ is a
  trace-class operator and satisfies
\begin{eqnarray*}
\Tr(\widetilde{P_\Lambda}P_\Lambda\,U_f) &=& (2\log(\Lambda) - B)_+(I(f))(1)\\
\Tr(\widetilde{P_\Lambda}P_\Lambda\,U_f) &=& 2\log(\Lambda)f(1) - H(f)(1) +
o(1)\\
\end{eqnarray*}
\end{thm}

For the last line we used that $B(I(f))(1) = H(I(f))(1) = H(f)(1)$ as $H =
\log(|x|) + \log(|y|)$ commutes with the Inversion $I$. The error is $o(1)$ for
$\Lambda \rightarrow\infty$ as it is bounded above in absolute value (assuming
$\Lambda > 1$) by
$$\sqrt{2\pi^2} \int_{|Y|\geq\Lambda^2} \log(|Y|)\;
\left|\Gamma(f)_a(Y)\right|\,dY$$ and $\Gamma(f)_a$ is a Schwartz function of
$Y\in\HH$. We note that if needed the Lemma {\bf\ref{L5}} gives to the term
$H(f)(1)$ a form more closely akin to the Weil's explicit formulae of number
theory. We note that Connes's computation in \cite{Co} also goes through an
intermediate stage essentially identical with  {\bf (C)} and that the
identification of the constant term with Weil's expression for the explicit
formula of number theory then requires a further discussion. The main result of
\cite{Bu1} and of this paper is thus the direct connection between $H$ and the
logarithmic derivatives of the Tate Gamma functions involved in the explicit
formulae.

{\bf Acknowledgements}

I thank the SSAS (``Soci\'et\'e de secours des amis des sciences'', quai de
Conti, Paris) for its financial support while this work was completed.

\vfill
\parskip = 0pt\parindent =200bp Jean-Fran\c cois Burnol\par
Universit\'e de Nice\--\ Sophia Antipolis\par Laboratoire J.-A. Dieudonn\'e\par
Parc Valrose\par F-06108 Nice C\'edex 02\par France\par burnol@math.unice.fr\par

\clearpage
\end{document}